\PassOptionsToPackage{unicode}{hyperref}
\PassOptionsToPackage{hyphens}{url}
\PassOptionsToPackage{dvipsnames,svgnames,x11names}{xcolor}
\documentclass[
  11pt,
]{article}
\usepackage{amsmath,amssymb}
\usepackage{iftex}
\ifPDFTeX
  \usepackage[T1]{fontenc}
  \usepackage[utf8]{inputenc}
  \usepackage{textcomp} 
\else 
  \usepackage{unicode-math} 
  \defaultfontfeatures{Scale=MatchLowercase}
  \defaultfontfeatures[\rmfamily]{Ligatures=TeX,Scale=1}
\fi
\usepackage{lmodern}
\ifPDFTeX\else
\fi
\IfFileExists{upquote.sty}{\usepackage{upquote}}{}
\IfFileExists{microtype.sty}{
  \usepackage[]{microtype}
  \UseMicrotypeSet[protrusion]{basicmath} 
}{}
\makeatletter
\@ifundefined{KOMAClassName}{
  \IfFileExists{parskip.sty}{%
    \usepackage{parskip}
  }{
    \setlength{\parindent}{0pt}
    \setlength{\parskip}{6pt plus 2pt minus 1pt}}
}{
  \KOMAoptions{parskip=half}}
\makeatother
\usepackage{xcolor}
\usepackage[margin=2.5cm]{geometry}
\usepackage{color}
\usepackage{fancyvrb}

\DefineVerbatimEnvironment{Highlighting}{Verbatim}{commandchars=\\\{\}}
\newenvironment{Shaded}{}{}

\newcommand{\BuiltInTok}[1]{\textcolor[rgb]{0.00,0.50,0.00}{#1}}

\newcommand{\CommentTok}[1]{\textcolor[rgb]{0.38,0.63,0.69}{\textit{#1}}}

\newcommand{\ControlFlowTok}[1]{\textcolor[rgb]{0.00,0.44,0.13}{\textbf{#1}}}

\newcommand{\DecValTok}[1]{\textcolor[rgb]{0.25,0.63,0.44}{#1}}

\newcommand{\ImportTok}[1]{\textcolor[rgb]{0.00,0.50,0.00}{\textbf{#1}}}

\newcommand{\KeywordTok}[1]{\textcolor[rgb]{0.00,0.44,0.13}{\textbf{#1}}}
\newcommand{\NormalTok}[1]{#1}
\newcommand{\OperatorTok}[1]{\textcolor[rgb]{0.40,0.40,0.40}{#1}}

\newcommand{\PreprocessorTok}[1]{\textcolor[rgb]{0.74,0.48,0.00}{#1}}

\newcommand{\SpecialCharTok}[1]{\textcolor[rgb]{0.25,0.44,0.63}{#1}}

\newcommand{\StringTok}[1]{\textcolor[rgb]{0.25,0.44,0.63}{#1}}

\usepackage{longtable,booktabs,array}
\usepackage{calc} 
\usepackage{etoolbox}
\makeatletter
\patchcmd\longtable{\par}{\if@noskipsec\mbox{}\fi\par}{}{}
\makeatother
\IfFileExists{footnotehyper.sty}{\usepackage{footnotehyper}}{\usepackage{footnote}}
\makesavenoteenv{longtable}
\usepackage{graphicx}
\makeatletter
\def\maxwidth{\ifdim\Gin@nat@width>\linewidth\linewidth\else\Gin@nat@width\fi}
\def\maxheight{\ifdim\Gin@nat@height>\textheight\textheight\else\Gin@nat@height\fi}
\makeatother
\setkeys{Gin}{width=\maxwidth,height=\maxheight,keepaspectratio}
\makeatletter
\def\fps@figure{htbp}
\makeatother
\providecommand{\tightlist}{%
  \setlength{\itemsep}{0pt}\setlength{\parskip}{0pt}}
\usepackage{amsmath,amssymb}
\ifLuaTeX
  \usepackage{selnolig}  
\fi
\IfFileExists{bookmark.sty}{\usepackage{bookmark}}{\usepackage{hyperref}}
\IfFileExists{xurl.sty}{\usepackage{xurl}}{} 
\hypersetup{
  pdftitle={Why \textbackslash bigl(4x(1-x)\textbackslash bigr)\^{}\{3/4\} approximates the binary entropy so well},
  pdfauthor={Jeroen van de Graaf},
  colorlinks=true,
  linkcolor={Maroon},
  filecolor={Maroon},
  citecolor={Blue},
  urlcolor={Blue},
  pdfcreator={LaTeX via pandoc}}

\title{Why \(\bigl(4x(1-x)\bigr)^{3/4}\) approximates the binary entropy
so well}
\usepackage{etoolbox}
\makeatletter
\providecommand{\subtitle}[1]{
  \apptocmd{\@title}{\par {\large #1 \par}}{}{}
}
\makeatother
\subtitle{A proof of Topsøe's convexity conjecture}
\author{Jeroen van de
Graaf\thanks{Departamento de Ciência da Computação, Universidade Federal de Minas Gerais (UFMG), Belo Horizonte, Brazil.}}
\date{November 2026}

\begin{document}
\maketitle
\begin{abstract}
Write the binary entropy as \(h(x) = u^{g(x)}\) with \(u = 4x(1-x)\).
Topsøe (2001) proved the sharp bounds \(1/\ln 4 \le g \le 1\),
i.e.~\(u \le h \le u^{1/\ln 4}\), and conjectured that the effective
exponent \(g\) is convex. We prove this conjecture. In the variable
\(s = (1-2x)^2\), \(g\) extends to a Pick function on
\(\mathbb{C} \setminus [1,\infty)\) with an explicit positive density on
the cut, so that \(g(x) = \sum_n \gamma_n (1-2x)^{2n}\) with
\(\gamma_n > 0\) and \((\gamma_n)_{n \ge 1}\) a moment sequence.
Pinelis's monotonicity theorem, Topsøe's power-type bounds and the
optimality of their exponents follow as corollaries, together with a
classification of all exponents \(\alpha\) for which \(u^\alpha\) bounds
\(h\). The proof is by hand except for one inequality in one real
variable on a bounded interval, certified by interval arithmetic. As an
application we give sharp, certified constants for the approximation
\(h \approx u^{3/4}\): \(-0.0065493 \le u^{3/4} - h \le 0.0080004\) on
\([0,1]\). Writing \(g = 1/\ln z\), with \(z\) between \(e\) and \(4\),
we also prove the one-sided variable-exponent bound
\(0 \le h - u^{1/\ln z} \le 0.0006377\) for
\(\ln z = 1 + (\ln 4 - 1)u^{1/10}\), about 170 times tighter than
\(h \ge u\).
\end{abstract}

The binary entropy \(h(x) = -x\log_2 x - (1-x)\log_2(1-x)\) can be
written exactly as a power of \(u = 4x(1-x)\) with a variable exponent,
\[
h(x) = u^{\,g(x)}, \qquad g(x) = \frac{\ln h(x)}{\ln\bigl(4x(1-x)\bigr)} .
\] Topsøe (2001) proved the sharp bounds \(u \le h \le u^{1/\ln 4}\),
i.e.~\(1/\ln 4 \le g \le 1\), and conjectured that \textbf{\(g\) is
convex} (his Conjecture 6.2). This note proves the conjecture.

\begin{itemize}
\tightlist
\item
  \textbf{Main result} (Theorem 1). In the variable \(s = (1-2x)^2\),
  \(g\) is a Pick function with an explicit positive density on the cut
  \([1,\infty)\). Hence \(g(x) = \sum_n \gamma_n (1-2x)^{2n}\) with
  \(\gamma_n > 0\), and \((\gamma_n)_{n\ge1}\) is a moment sequence.
\item
  \textbf{Consequences} (Corollary 7). \(g\) is strictly convex;
  Pinelis's monotonicity theorem, Topsøe's power bounds and the
  optimality of their exponents follow at once; for every exponent
  strictly between \(1/\ln 4\) and \(1\), \(u^\alpha\) crosses \(h\)
  exactly once on each side of \(\tfrac12\).
\item
  \textbf{An application} (Section 5). The constant exponent \(3/4\)
  gives the approximation \(h \approx \bigl(4x(1-x)\bigr)^{3/4}\) with
  \[
  -0.0065493 \;\le\; \bigl(4x(1-x)\bigr)^{3/4} - h(x) \;\le\; 0.0080004 \qquad (0 \le x \le 1).
  \]
\end{itemize}

The proof of the main result is by hand except for one inequality in one
real variable on a bounded interval (Lemma 5, middle range), which has a
margin of \(0.65\) and is certified by interval arithmetic. The sharp
constants for the \(3/4\) approximation are also certified by interval
arithmetic.

\hypertarget{setting-and-main-result}{%
\section{Setting and main result}\label{setting-and-main-result}}

Throughout, \(x \in [0,1]\) and \[
u = 4x(1-x), \qquad s = (1-2x)^2 = 1-u, \qquad c = \frac{1}{2\ln 2} = \frac{1}{\ln 4} \approx 0.72135 .
\] The point \(x = \tfrac12\) corresponds to \(s = 0\), and
\(x \in \{0,1\}\) to \(s = 1\). Everything is symmetric under
\(x \mapsto 1-x\). The Taylor series of \(h\) at \(\tfrac12\) is \[
h = 1 - c \sum_{n\ge1} \frac{s^n}{n(2n-1)}, \qquad \sum_{n\ge1}\frac{1}{n(2n-1)} = \frac1c , \tag{1}
\] the identity coming from \(h = 0\) at \(s = 1\). In the variable
\(s\) the effective exponent is \[
g(x) = G(s), \qquad G(s) = \frac{\ln h}{\ln(1-s)},
\] with \(g(\tfrac12) = G(0) = c\) and \(g(0) = g(1) = 1\) by
continuity. Write \(\Omega = \mathbb{C} \setminus [1,\infty)\),
\(\mathbb{H} = \{\operatorname{Im} s > 0\}\), and \(\operatorname{Ln}\)
for the principal logarithm.

\textbf{Theorem 1.} \(G\) extends analytically to \(\Omega\), and \[
G(s) = \frac12 + \frac1\pi \int_1^\infty \frac{w(\sigma)}{\sigma - s}\, d\sigma \qquad (s \in \Omega),
\] where \(w\) is continuous and strictly positive on \((1,\infty)\),
given explicitly in Lemma 4, with \(w(\sigma) = O(1/\ln^2\sigma)\) as
\(\sigma \to \infty\) and \(w(\sigma) \to 0\) as
\(\sigma \downarrow 1\). Consequently, uniformly for \(x \in [0,1]\), \[
g(x) = \sum_{n \ge 0} \gamma_n\, (1-2x)^{2n}, \qquad \gamma_0 = \frac{1}{\ln 4}, \qquad \gamma_n = \frac1\pi \int_1^\infty w(\sigma)\, \sigma^{-n-1}\, d\sigma > 0 \quad (n \ge 1),
\] and \(\sum_{n\ge0} \gamma_n = 1\).

In particular \(G\) is a Pick (Nevanlinna) function:
\(\operatorname{Im} G > 0\) on \(\mathbb{H}\). The consequences,
including Topsøe's conjecture, are collected in Corollary 7.

\hypertarget{prior-work}{%
\section{Prior work}\label{prior-work}}

\textbf{Topsøe's bounds.} Topsøe {[}T01{]} proved two pairs of bounds
(stated there in nats; we convert to bits):

\begin{longtable}[]{@{}
  >{\raggedright\arraybackslash}p{(\columnwidth - 4\tabcolsep) * \real{0.3333}}
  >{\raggedright\arraybackslash}p{(\columnwidth - 4\tabcolsep) * \real{0.3333}}
  >{\raggedright\arraybackslash}p{(\columnwidth - 4\tabcolsep) * \real{0.3333}}@{}}
\toprule\noalign{}
\begin{minipage}[b]{\linewidth}\raggedright
Name
\end{minipage} & \begin{minipage}[b]{\linewidth}\raggedright
Bound
\end{minipage} & \begin{minipage}[b]{\linewidth}\raggedright
Sharpness
\end{minipage} \\
\midrule\noalign{}
\endhead
\bottomrule\noalign{}
\endlastfoot
Logarithmic &
\(\ln 2 \cdot \log_2 x \log_2(1-x) \le h \le \log_2 x \log_2(1-x)\) &
best constants \\
Power-type & \(4x(1-x) \le h \le \bigl(4x(1-x)\bigr)^{1/\ln 4}\) & best
exponents \\
\end{longtable}

According to {[}T01{]}, the lower power bound is a special case of
Theorem 2.6 of Harremoës and Topsøe {[}HT01{]}, where the power bounds
and their generalizations to more than two points are discussed, and
\(1/\ln 4\) is the largest exponent possible in the upper bound. The
upper bound improves the exponent \(\tfrac12\) of Lin {[}L91{]}. In
{[}T01{]}:

\begin{itemize}
\tightlist
\item
  the lower power bound is proved from the convexity of \(\ln x/(x-1)\);
\item
  the upper power bound is proved by expanding the difference
  \(u^{c} - h\) in powers of \((1-2x)^2\). Topsøe notes that the series
  has terms of both signs and so cannot be used directly, and concludes
  by showing that the second derivative vanishes only once;
\item
  the printed proof of the logarithmic bounds is due to Arjomand,
  Bahramgiri and Rouhani; a conjecture on related Taylor coefficients
  (Conjecture 6.1) was settled independently by Berg {[}B01{]}, using
  Pick functions, and by Laczkovich.
\end{itemize}

\textbf{The effective exponent.} Topsøe {[}T01, (6.6){]} introduces
\(g\), observes that the power bounds are equivalent to
\(g(\tfrac12) \le g \le g(0)\), calls his proof ``somewhat ad hoc'', and
states Conjecture 6.2: \emph{the function \(g\) is convex}. He reports
that Pinelis proved \(g\) monotone on \([0,\tfrac12]\) and on
\([\tfrac12,1]\) with an l'Hospital-type rule for monotonicity; see also
{[}P02{]}. Pinelis {[}P01, §2.2{]} returns to the convexity conjecture:
he reduces it to an algorithmic analysis of the oscillation pattern of a
rational function, but does not carry it out, judging the calculations
too large. We are not aware of a later resolution.

\textbf{The exponent 3/4.} The author found the approximation
\(h \approx u^{3/4}\) around 1995 in numerical experiments with
Mathematica, during the work with C. A. Fuchs that led to {[}FG99{]}.
That paper uses the bound \(h \le 2\sqrt{x(1-x)} = u^{1/2}\) to relate
distinguishability measures, plots it against \(h\) {[}FG99, Fig. 1{]},
and remarks that the bounding function resembles \(h\) quite well. The
exponent \(3/4\) arose from looking for a better exponent than
\(\tfrac12\); the observation was not included in {[}FG99{]} and was
never published. The approximation has since appeared informally,
without proof, in two places. In 2015 J. D'Aurizio remarked on
Mathematics Stack Exchange {[}D15{]} that \(h \approx u^{3/4}\) is
neither an upper nor a lower bound and that the absolute error is always
below \(8\cdot10^{-3}\); by Theorem 8 the maximum is in fact
\(0.0080003\), attained near \(x \approx 0.0089\). A 2021 question there
{[}S21{]} asked for a proof of \(h \le u^{3/4} + 0.007\); the answer it
received points to Topsøe's bounds instead and gives no proof. In 2019 a
blog post by J. D. Cook {[}C19{]}, crediting a reader (S. Visscher),
used the approximation to invert the binary entropy. We are not aware of
a published analysis of it.

\hypertarget{proof-of-theorem-1}{%
\section{Proof of Theorem 1}\label{proof-of-theorem-1}}

\hypertarget{analytic-continuation-of-h}{%
\subsection{\texorpdfstring{Analytic continuation of
\(h\)}{Analytic continuation of h}}\label{analytic-continuation-of-h}}

\textbf{Lemma 2.} Put \(q(s) = s \int_0^1 \frac{d\tau}{1 - s\tau^2}\),
so that
\(q(s) = \sqrt{s}\,\operatorname{artanh}\sqrt{s} = \sum_{n\ge1} \frac{s^n}{2n-1}\)
for \(|s| < 1\). Then for \(s \in [0,1)\) \[
h = 1 + c \int_0^1 \frac{\operatorname{Ln}(1 - s\tau^2)}{\tau^2}\, d\tau = 1 - c\, \operatorname{Ln}(1-s) - 2c\, q(s). \tag{2}
\] Both right-hand sides are analytic on \(\Omega\) and define the
continuation of \(h\); it satisfies \(h(\bar s) = \overline{h(s)}\).

\emph{Proof.} Differentiating (1),
\(h'(s) = -c\sum_{n\ge1} \frac{s^{n-1}}{2n-1} = -c \int_0^1 \frac{d\tau}{1 - s\tau^2}\).
Integrating from \(h(0) = 1\) gives the integral form, and integrating
by parts in \(\tau\) gives the closed form. (With
\(t = \sqrt s = 1 - 2x\), the closed form is the familiar
\(\ln 2 \cdot h = \ln 2 - \tfrac12\ln(1-t^2) - t\operatorname{artanh} t\).)
For \(s \in \Omega\) and \(\tau \in (0,1]\),
\(1 - s\tau^2 \notin (-\infty, 0]\), so the integrands are analytic in
\(s\). \(\square\)

\textbf{Lemma 3.} \(\operatorname{Im} h(s) < 0\) for
\(s \in \mathbb{H}\), and \(h(s) > 0\) for \(s \in (-\infty, 1)\). Hence
\(h\) has no zeros in \(\Omega\), \(\operatorname{Ln} h\) is analytic on
\(\Omega\), and \(G = \operatorname{Ln} h / \operatorname{Ln}(1-s)\) is
analytic on \(\Omega\), real on \((-\infty,1)\), with
\(G(\bar s) = \overline{G(s)}\).

\emph{Proof.} For \(s \in \mathbb{H}\),
\(\operatorname{Im}(1 - s\tau^2) = -\tau^2 \operatorname{Im} s < 0\), so
every logarithm in the integral form of (2) has imaginary part in
\((-\pi, 0)\). On \((-\infty,1)\),
\(h' = -c\int_0^1 d\tau/(1-s\tau^2) < 0\) and \(h(1^-) = 0\), so
\(h > 0\). Thus \(h(\Omega)\) avoids \((-\infty,0]\). The denominator
\(\operatorname{Ln}(1-s)\) vanishes on \(\Omega\) only at \(s = 0\),
where the numerator vanishes too (\(h(0) = 1\)), so the singularity is
removable. \(\square\)

\hypertarget{boundary-values-on-the-cut}{%
\subsection{Boundary values on the
cut}\label{boundary-values-on-the-cut}}

\textbf{Lemma 4.} For \(\sigma > 1\) write \(\sigma = (1+2y)^2\) with
\(y > 0\), and put \[
R = y\ln(1+\tfrac1y) + \ln(1+y) > 0, \qquad v = 4y(1+y) = \sigma - 1, \qquad \varphi = \arctan\frac{\pi y}{R} \in (0, \tfrac\pi2).
\] As \(s \to \sigma\) from \(\mathbb{H}\), \(h\) and
\(\operatorname{Ln}(1-s)\) have the limits \[
\ln 2 \cdot h(\sigma + i0) = -R - i\pi y, \qquad \operatorname{Ln}(1 - \sigma - i0) = \ln v - i\pi,
\] so \(\operatorname{Ln} h(\sigma + i0) = \ln|h| - i(\pi - \varphi)\)
with \(|h| = \sqrt{R^2 + \pi^2 y^2}/\ln 2\), and \(G\) extends
continuously from \(\mathbb{H}\) to \((1,\infty)\) with \[
w(\sigma) := \operatorname{Im} G(\sigma + i0) = \frac{\Phi(y)}{\ln^2 v + \pi^2}, \qquad \Phi(y) = \pi \ln\frac{|h|}{v} + \varphi \ln v .
\]

\emph{Proof.} As \(s \to \sigma\) from \(\mathbb{H}\),
\(t = \sqrt s \to T = 1 + 2y\) with \(\operatorname{Im} t > 0\). Then
\(\frac{1+t}{1-t}\) approaches the negative real number
\(-\frac{T+1}{T-1}\) from above, so
\(\operatorname{artanh} t = \tfrac12\operatorname{Ln}\frac{1+t}{1-t} \to \tfrac12 \ln\frac{T+1}{T-1} + \tfrac{i\pi}{2}\),
and \(\operatorname{Ln}(1-s) \to \ln(\sigma - 1) - i\pi\). Substituting
in the closed form of (2), \[
h(\sigma+i0) = 1 - c\ln v - cT\ln\tfrac{T+1}{T-1} - i\pi c\,(T-1),
\] and with \(1 = 2c\ln 2\), \(T - 1 = 2y\),
\(\frac{T+1}{T-1} = \frac{1+y}{y}\) the real part simplifies to
\(-R/\ln 2\) and the imaginary part to \(-\pi y/\ln 2\). Since
\(-R - i\pi y\) lies in the open third quadrant, its argument is
\(-\pi + \varphi\). Finally
\(\operatorname{Im}\frac{A + i\theta}{a - i\pi} = \frac{\theta a + \pi A}{a^2 + \pi^2}\)
with \(A = \ln|h|\), \(\theta = \varphi - \pi\), \(a = \ln v\).
Continuity up to the cut holds because the limits of \(h\) avoid
\((-\infty, 0]\) and those of \(\operatorname{Ln}(1-s)\) avoid \(0\).
\(\square\)

\textbf{Lemma 5.} \(\Phi(y) > 0\) for all \(y > 0\). (Numerically
\(\min \Phi \approx 0.6531\), at \(y \approx 0.168\).)

\emph{Proof.} Three ranges.

\begin{itemize}
\tightlist
\item
  \emph{Small \(y\), \(0 < y \le 10^{-4}\).} Let
  \(L = \ln(1/y) \ge 9.21\). Then \(R \ge yL\), so \(|h| \ge yL/\ln 2\),
  while \(v \le 4y(1 + 10^{-4})\). Also
  \(\varphi \le \pi y/R \le \pi/L\) and \(0 < -\ln v \le L\), so
  \(\varphi \ln v \ge -\pi\). Hence \[
  \Phi \ge \pi\Bigl[\ln \frac{L}{4\ln 2\,(1+10^{-4})} - 1\Bigr] > 0,
  \] because \(L \ge 9.21 > 4e\ln 2\,(1+10^{-4}) \approx 7.54\).
\item
  \emph{Large \(y\), \(y \ge 10\).} With
  \(\psi = \frac\pi2 - \varphi = \arctan\frac{R}{\pi y}\), one has
  \(\Phi = \frac\pi2 \ln\frac{|h|^2}{v} - \psi \ln v\). Since
  \(\ln^2 2\,|h|^2 \ge \pi^2y^2\) and
  \(\frac{y}{1+y} \ge \frac{10}{11}\), the first term is at least
  \(\frac\pi2\ln\bigl(\frac{\pi^2}{4\ln^2 2}\cdot\frac{10}{11}\bigr) > 2.42\).
  For the second, put \(\lambda = \ln(1+y)\). Since
  \(y\ln(1+1/y) \le 1\), \(R \le 1 + \lambda\); also
  \(0 < \ln v \le \ln 4 + 2\lambda\) and \(\psi \le R/(\pi y)\), so \[
  \psi \ln v \;\le\; B(y) := \frac{(1+\lambda)(\ln 4 + 2\lambda)}{\pi y} .
  \] \(B\) is decreasing on \([10,\infty)\): its logarithmic derivative
  is
  \(\bigl[\frac{1}{1+\lambda} + \frac{2}{\ln 4 + 2\lambda}\bigr]\frac{1}{1+y} - \frac1y\),
  which is negative as soon as the bracket is at most \(1\),
  i.e.~\(2\lambda^2 - (2 - \ln 4)\lambda - 2 \ge 0\); this holds for
  \(\lambda \ge \ln 11 \approx 2.398\). Since \(B(10) \approx 0.6687\),
  \(\Phi > 2.42 - 0.67 = 1.75\).
\item
  \emph{Middle range, \(10^{-4} \le y \le 10\).} Certified by interval
  arithmetic: Arb ball arithmetic {[}Arb{]} with adaptive bisection
  shows \(\Phi > 0\) on 92 subintervals covering the range (Appendix A).
  \(\square\)
\end{itemize}

\begin{figure}
\centering
\includegraphics[width=0.8\textwidth,height=\textheight]{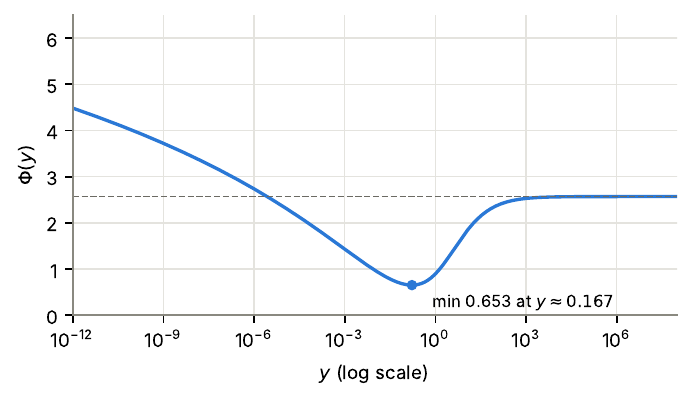}
\caption{The function \(\Phi(y)\) on the cut, log scale in \(y\). It
stays above \(0.65\) and tends to
\(\pi\ln\frac{\pi}{2\ln 2} \approx 2.570\) as \(y \to \infty\).}
\end{figure}

\hypertarget{behaviour-at-s-1-and-at-infinity}{%
\subsection{\texorpdfstring{Behaviour at \(s = 1\) and at
infinity}{Behaviour at s = 1 and at infinity}}\label{behaviour-at-s-1-and-at-infinity}}

\textbf{Lemma 6.} \(G(s) \to 1\) as \(s \to 1\) and
\(G(s) \to \tfrac12\) as \(|s| \to \infty\), in both cases for
\(s \in \Omega\) and uniformly in the direction of approach.
Consequently \(G\) is bounded on \(\mathbb{H}\).

\emph{Proof.} \emph{Near \(s = 1\).} Put \(\varepsilon = 1 - \sqrt s\).
For \(s \in \Omega\) near \(1\), \(\varepsilon \to 0\) and
\(\varepsilon \notin (-\infty, 0]\), so
\(\arg\varepsilon \in (-\pi,\pi)\); since \(2 - \varepsilon\) has a
small argument of the sign opposite to that of \(\varepsilon\),
\(\operatorname{Ln}(1-s) = \operatorname{Ln}\varepsilon + \operatorname{Ln}(2-\varepsilon)\),
and
\(\operatorname{artanh}\sqrt s = \tfrac12[\operatorname{Ln}(2-\varepsilon) - \operatorname{Ln}\varepsilon]\).
The closed form of (2) becomes, exactly, \[
h = 1 - c\,\varepsilon \operatorname{Ln}\varepsilon - c\,(2-\varepsilon)\operatorname{Ln}(2-\varepsilon) = c\,\varepsilon\,\bigl(\Lambda + 1 + \ln 2 + O(\varepsilon)\bigr), \qquad \Lambda = \operatorname{Ln}\tfrac1\varepsilon ,
\] using \(1 = 2c\ln 2\). Here
\(|\Lambda| \ge \ln\frac{1}{|\varepsilon|} \to \infty\) uniformly in
\(\arg\varepsilon\), so the last factor is nonzero for small
\(|\varepsilon|\) and its logarithm is \(O(\ln|\Lambda|)\). The
principal logarithm of a product of three factors equals the sum of
their principal logarithms up to an integer multiple of \(2\pi i\) of
absolute value at most \(2\pi\), so \[
\operatorname{Ln} h = \operatorname{Ln}\varepsilon + O(\ln|\Lambda|), \qquad \operatorname{Ln}(1-s) = \operatorname{Ln}\varepsilon + O(1),
\] uniformly in \(\arg \varepsilon\). As
\(|\operatorname{Ln}\varepsilon| = |\Lambda| \to \infty\),
\(G = 1 + O(\ln|\Lambda|/|\Lambda|) \to 1\).

\emph{Near infinity.} Let \(t = \sqrt s\), so
\(\operatorname{Re} t \ge 0\), \(|t| \to \infty\), and
\(\operatorname{Im} t \ne 0\) as soon as \(|t| > 1\). For
\(\pm\operatorname{Im} t > 0\),
\(\operatorname{artanh} t = \operatorname{artanh}\frac1t \pm \frac{i\pi}{2}\),
hence
\(t\operatorname{artanh} t = \pm\frac{i\pi}{2}t + 1 + O(|t|^{-2})\) and,
by (2), \[
h = \mp\, i\pi c\,\sqrt s - c\operatorname{Ln}(1-s) + O(1) = \mp\, i\pi c\,\sqrt s\,\bigl(1 + O(|s|^{-1/2}\ln|s|)\bigr),
\] uniformly in \(\arg s\). Thus
\(\operatorname{Ln} h = \tfrac12 \operatorname{Ln} s + O(1)\) and
\(\operatorname{Ln}(1-s) = \operatorname{Ln} s + O(1)\), so
\(G \to \tfrac12\).

\emph{Boundedness.} \(G\) is analytic across \((-\infty,1)\) and extends
continuously from \(\mathbb{H}\) to \((1,\infty)\) (Lemma 4). With the
two limits just proved, \(G\) extends continuously to the closure of
\(\mathbb{H}\) in the Riemann sphere, which is compact. \(\square\)

\hypertarget{conclusion-of-the-proof}{%
\subsection{Conclusion of the proof}\label{conclusion-of-the-proof}}

By Lemmas 4--6, \(w\) is continuous and positive on \((1,\infty)\), and
bounded: as \(\sigma\to\infty\), \(\Phi\) stays bounded, so
\(w = O(1/\ln^2\sigma)\); as \(\sigma \downarrow 1\),
\(w \to \operatorname{Im} 1 = 0\). Hence \[
F(s) = \frac1\pi\int_1^\infty \frac{w(\sigma)}{\sigma - s}\, d\sigma
\] converges absolutely for \(s \in \Omega\) (note
\(\int^\infty \frac{d\sigma}{\sigma\ln^2\sigma} < \infty\)), is analytic
there, is real on \((-\infty,1)\), and tends to \(0\) as
\(s \to -\infty\) by dominated convergence. On \(\mathbb{H}\), \[
\operatorname{Im} F(s) = \frac1\pi\int_1^\infty w(\sigma)\, \frac{\operatorname{Im} s}{|\sigma - s|^2}\, d\sigma
\] is the Poisson integral of the bounded function \(\tilde w\) equal to
\(w\) on \((1,\infty)\) and to \(0\) on \((-\infty,1]\). Two standard
facts about the Poisson integral of a bounded function {[}Ga81, Ch. I,
§3{]} apply: it is bounded by \(\sup|\tilde w|\), and it tends to
\(\tilde w(\xi)\) at every point \(\xi\) where \(\tilde w\) is
continuous, here every \(\xi \in \mathbb{R}\setminus\{1\}\).

Now let \(U = \operatorname{Im} G - \operatorname{Im} F\) on
\(\mathbb{H}\). It is harmonic and bounded (Lemma 6), and it tends to
\(0\) at every boundary point \(\xi \in \mathbb{R}\setminus\{1\}\): on
\((1,\infty)\) both terms tend to \(w\) (Lemma 4), and on
\((-\infty,1)\) both tend to \(0\) (Lemma 3). By the Phragmén--Lindelöf
form of the maximum principle (a bounded harmonic function whose
boundary limits are \(\le 0\) except at finitely many boundary points,
the point at infinity included, is \(\le 0\); see e.g.~{[}Ra95, Thm
3.6.9{]}), applied to \(U\) and \(-U\) with exceptional points \(1\) and
\(\infty\), \(U \equiv 0\).

So \(G - F\) is analytic on \(\mathbb{H}\) with vanishing imaginary
part, hence equal to a real constant \(\alpha\) there, and by the
identity theorem on all of the connected set \(\Omega\). Letting
\(s \to -\infty\) gives \(\alpha = \tfrac12\) (Lemma 6). This proves the
representation, and \(\operatorname{Im} G = \operatorname{Im} F > 0\) on
\(\mathbb{H}\).

For \(|s| \le r < 1\), expand
\(\frac{1}{\sigma - s} = \sum_{n\ge0} s^n \sigma^{-n-1}\), which
converges uniformly in \(\sigma \ge 1\); integrating term by term gives
the Taylor coefficients of \(G\) at \(0\):
\(\gamma_n = \frac1\pi\int_1^\infty w\,\sigma^{-n-1}d\sigma > 0\) for
\(n \ge 1\), and \(\gamma_0 = G(0) = c\). Since all \(\gamma_n\) are
positive and \(G(s) \to 1\) as \(s \uparrow 1\), monotone convergence
gives \(\sum_n \gamma_n = 1\). The series \(\sum\gamma_n s^n\) therefore
converges uniformly on \([0,1]\) by the Weierstrass test, and equals
\(g(x)\) with \(s = (1-2x)^2\). \(\square\)

\hypertarget{consequences}{%
\section{Consequences}\label{consequences}}

\textbf{Corollary 7.}

\begin{enumerate}
\def\labelenumi{\arabic{enumi}.}
\tightlist
\item
  \textbf{(Topsøe's Conjecture 6.2.)} \(g\) is strictly convex on
  \([0,1]\).
\item
  \textbf{(Pinelis's theorem.)} \(g\) is strictly decreasing on
  \([0,\tfrac12]\) and strictly increasing on \([\tfrac12,1]\).
\item
  \textbf{(Topsøe's power bounds.)} \(u \le h \le u^{1/\ln 4}\), with
  equality only at \(x \in \{0,\tfrac12,1\}\).
\item
  \textbf{(Classification of exponents.)} Let \(\alpha > 0\). Then
  \(u^\alpha \ge h\) on \([0,1]\) if and only if \(\alpha \le 1/\ln 4\);
  \(u^\alpha \le h\) on \([0,1]\) if and only if \(\alpha \ge 1\); and
  for \(1/\ln 4 < \alpha < 1\) the curves \(u^\alpha\) and \(h\) cross
  exactly once in \((0,\tfrac12)\), at a point \(x^\star(\alpha)\), with
  \(h < u^\alpha\) on \((0,x^\star)\) and \(h > u^\alpha\) on
  \((x^\star,\tfrac12)\).
\item
  \textbf{(Curvature at the centre.)}
  \(g''(\tfrac12) = 8\gamma_1 = 4c^2 - \tfrac83 c \approx 0.1578\).
\item
  \textbf{(Moments.)} With \(\tau = 1/\sigma\),
  \(\gamma_n = \frac1\pi\int_0^1 \tau^{n-1}\, w(1/\tau)\, d\tau\), so
  \((\gamma_n)_{n\ge1}\) is a Hausdorff moment sequence and in
  particular completely monotone.
\end{enumerate}

\emph{Proof.} Each function \((1-2x)^{2n}\), \(n \ge 1\), is strictly
convex on \([0,1]\) and strictly decreasing on \([0,\tfrac12]\). A
uniformly convergent sum of such functions with positive coefficients,
plus a constant, has both properties. This gives 1 and 2.

For 3 and 4, note that for \(0 < u < 1\),
i.e.~\(x \notin \{0,\tfrac12,1\}\), we have \(\ln u < 0\) and
\(h = u^{g}\), so \(h \le u^\alpha\) if and only if \(g(x) \ge \alpha\),
with the same equivalence for strict inequalities and for the reversed
ones. By 2, \(g\) maps \((0,\tfrac12)\) continuously, strictly
monotonically and onto \((c, 1)\). So \(g > c\) and \(g < 1\) there,
which is 3; \(g \ge \alpha\) on all of \((0,\tfrac12)\) if and only if
\(\alpha \le c\), and \(g \le \alpha\) if and only if \(\alpha \ge 1\);
and for \(c < \alpha < 1\) the equation \(g = \alpha\) has exactly one
solution \(x^\star\), with \(g > \alpha\) to its left. This is 4.

For 5, only the \(n = 1\) term contributes to \(g''(\frac12)\), and
expanding \(G\) at \(s = 0\) with (1) gives
\(\gamma_1 = c(\frac c2 - \frac13)\). Item 6 is the substitution; the
measure \(w(1/\tau)\,d\tau/\pi\) is finite because
\(\int_1^\infty w\,\sigma^{-2}d\sigma < \infty\). \(\square\)

\begin{figure}
\centering
\includegraphics[width=0.8\textwidth,height=\textheight]{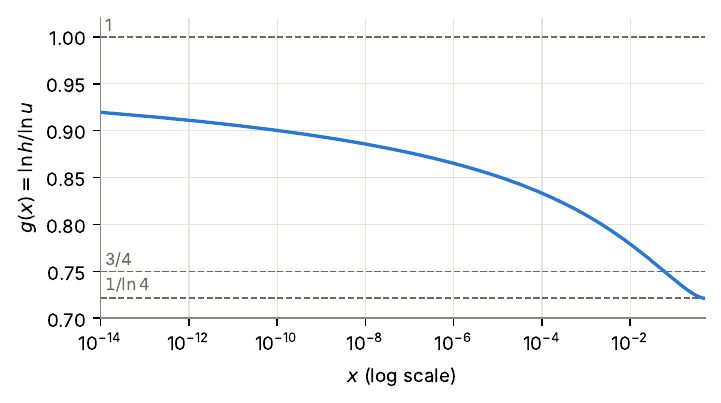}
\caption{Effective exponent \(g(x) = \ln h / \ln u\) (log scale in
\(x\)). It is convex as a function of \(x\); on this logarithmic axis
that is not visible.}
\end{figure}

The exponent \(g\) approaches \(1\) very slowly near the endpoints:

\begin{longtable}[]{@{}
  >{\raggedright\arraybackslash}p{(\columnwidth - 14\tabcolsep) * \real{0.0968}}
  >{\centering\arraybackslash}p{(\columnwidth - 14\tabcolsep) * \real{0.1290}}
  >{\centering\arraybackslash}p{(\columnwidth - 14\tabcolsep) * \real{0.1290}}
  >{\centering\arraybackslash}p{(\columnwidth - 14\tabcolsep) * \real{0.1290}}
  >{\centering\arraybackslash}p{(\columnwidth - 14\tabcolsep) * \real{0.1290}}
  >{\centering\arraybackslash}p{(\columnwidth - 14\tabcolsep) * \real{0.1290}}
  >{\centering\arraybackslash}p{(\columnwidth - 14\tabcolsep) * \real{0.1290}}
  >{\centering\arraybackslash}p{(\columnwidth - 14\tabcolsep) * \real{0.1290}}@{}}
\toprule\noalign{}
\begin{minipage}[b]{\linewidth}\raggedright
\(x\)
\end{minipage} & \begin{minipage}[b]{\linewidth}\centering
\(0.5\)
\end{minipage} & \begin{minipage}[b]{\linewidth}\centering
\(0.25\)
\end{minipage} & \begin{minipage}[b]{\linewidth}\centering
\(0.1\)
\end{minipage} & \begin{minipage}[b]{\linewidth}\centering
\(10^{-2}\)
\end{minipage} & \begin{minipage}[b]{\linewidth}\centering
\(10^{-3}\)
\end{minipage} & \begin{minipage}[b]{\linewidth}\centering
\(10^{-6}\)
\end{minipage} & \begin{minipage}[b]{\linewidth}\centering
\(10^{-12}\)
\end{minipage} \\
\midrule\noalign{}
\endhead
\bottomrule\noalign{}
\endlastfoot
\(g(x)\) & \(0.7213\) & \(0.7270\) & \(0.7411\) & \(0.7792\) &
\(0.8100\) & \(0.8652\) & \(0.9110\) \\
\end{longtable}

For instance \(u^{0.9} \le h\) fails only for \(x \lesssim 10^{-10}\),
and \(u^{0.95} \le h\) only for \(x \lesssim 10^{-28}\).

\textbf{Remarks.}

\begin{enumerate}
\def\labelenumi{(\roman{enumi})}
\item
  The method is the one Berg used for Topsøe's Conjecture 6.1 {[}B01{]}:
  show that a function related to the entropy is a Pick function and
  read off properties of its Taylor coefficients from the integral
  representation.
\item
  The only computer-assisted step is the middle range of Lemma 5, an
  inequality in one real variable with margin \(0.65\). It could be
  replaced by a hand argument on a few subintervals.
\item
  Numerical checks, independent of the proof: the first 398 coefficients
  \(\gamma_n\) computed from the Taylor series in 60-digit arithmetic
  are all positive, satisfy \((-1)^k\Delta^k\gamma_n \ge 0\) for
  \(k \le 15\), and agree to 10 digits with the integral formula;
  \(\operatorname{Im} G > 0\) at some 3000 points of \(\mathbb{H}\) with
  \(|s|\) and \(|s-1|\) ranging from \(10^{-30}\) to \(10^{30}\); and a
  direct evaluation gives \(g'' \ge 0.1577\) on a grid in \(x\) reaching
  \(10^{-20}\).
\end{enumerate}

\hypertarget{the-approximation-h-approx-u34}{%
\section{\texorpdfstring{The approximation
\(h \approx u^{3/4}\)}{The approximation h \textbackslash approx u\^{}\{3/4\}}}\label{the-approximation-h-approx-u34}}

Replacing the slowly varying exponent \(g(x)\) by the constant \(3/4\)
gives a simple and accurate approximation. By Corollary 7(4),
\(u^{3/4}\) crosses \(h\) exactly once on each side of \(\tfrac12\), at
\(x^\star \approx 0.059277\) and \(1 - x^\star\): it overestimates \(h\)
near the endpoints and underestimates it in the middle.

\begin{figure}
\centering
\includegraphics[width=0.75\textwidth,height=\textheight]{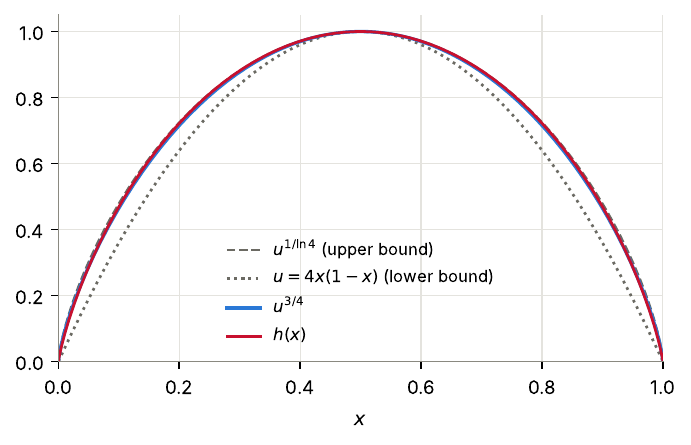}
\caption{The binary entropy \(h\) (red), the approximation \(u^{3/4}\)
(blue) and Topsøe's bounds (gray). At this scale \(h\) and \(u^{3/4}\)
are indistinguishable.}
\end{figure}

\begin{figure}
\centering
\includegraphics[width=0.75\textwidth,height=\textheight]{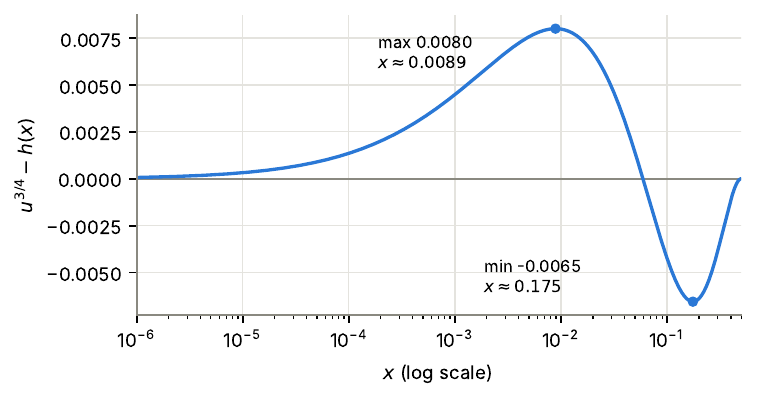}
\caption{Error \(u^{3/4} - h\) on \((0, \tfrac12]\), log scale in \(x\).
The largest error is a narrow spike close to the endpoint.}
\end{figure}

\hypertarget{a-sharp-computer-assisted-bound}{%
\subsection{A sharp, computer-assisted
bound}\label{a-sharp-computer-assisted-bound}}

\textbf{Theorem 8.} For all \(x \in [0,1]\), \[
-0.0065493 \;\le\; \bigl(4x(1-x)\bigr)^{3/4} - h(x) \;\le\; 0.0080004 .
\] The extreme values are \(\approx -0.0065492\) at
\(x \approx 0.175235\) and \(\approx +0.0080003\) at
\(x \approx 0.008892\) (and their mirror images), so the constants are
sharp to about \(10^{-7}\). In particular \(h \le u^{3/4} + 0.0065493\),
which answers the question in {[}S21{]}, and
\(|u^{3/4} - h| \le 0.0080004\), in agreement with the value
\(8\cdot10^{-3}\) given in {[}D15{]}.

\emph{Proof.} Write \(D = u^{3/4} - h\). By symmetry take
\(x \in [0,\tfrac12]\).

\begin{itemize}
\tightlist
\item
  On \([0, 10^{-4}]\), both \(u^{3/4}\) and \(h\) are increasing and
  nonnegative, so \(-h(10^{-4}) \le D \le (4\cdot 10^{-4})^{3/4}\),
  i.e.~\(D \in [-0.0015, 0.0029]\).
\item
  On \([10^{-4}, \tfrac12]\), \(D\) is \(C^1\). The interval is split
  adaptively into subintervals \(X\) with midpoint \(m\), and on each
  the mean-value form \(D(X) \subseteq D(m) + D'(X)\,(X - m)\) is
  evaluated in Arb ball arithmetic {[}Arb{]} with outward rounding. A
  subinterval is accepted only if the whole enclosure lies strictly
  inside the bounds; 145 subintervals suffice (Appendix B). \(\square\)
\end{itemize}

\textbf{Remark (binomial coefficients).} For integers \(0 < k < n\) and
\(p = k/n\), the classical bounds {[}CT06, Ch. 11; MS77, Ch. 10, Lemma
7{]} \[
\frac{2^{n h(p)}}{n+1} \le \binom{n}{k} \le 2^{n h(p)}, \qquad \frac{2^{n h(p)}}{\sqrt{8np(1-p)}} \le \binom{n}{k} \le \frac{2^{n h(p)}}{\sqrt{2\pi np(1-p)}}
\] combine with Corollary 7(3) and Theorem 8, because
\(t \mapsto 2^{nt}\) is increasing. With \(u = 4p(1-p)\), so that
\(8np(1-p) = 2nu\) and \(2\pi np(1-p) = \pi n u/2\), this gives \[
\frac{2^{nu}}{n+1} \le \binom{n}{k} \le 2^{n u^{1/\ln 4}}, \qquad \frac{2^{n(u^{3/4} - 0.0080004)}}{\sqrt{2nu}} \le \binom{n}{k} \le \frac{2^{n(u^{3/4} + 0.0065493)}}{\sqrt{\pi n u/2}} .
\] These are weaker than the bounds in terms of \(h\), but involve only
elementary functions of \(k/n\). In the second pair the exponent is
accurate to within \(0.008\,n\) bits, which suffices for growth rates
but not for the exact size of \(\binom{n}{k}\) when \(n\) is large.

\hypertarget{an-elementary-bound}{%
\subsection{An elementary bound}\label{an-elementary-bound}}

A weaker bound needs no computer. For \(0 < \alpha < 1\) the binomial
series gives \[
u^\alpha = (1-s)^\alpha = 1 - \sum_{n\ge1} b_n s^n, \qquad b_1 = \alpha,\quad b_{n+1} = b_n\,\frac{n-\alpha}{n+1},
\] with \(b_n > 0\) and
\(b_n \sim \frac{\alpha}{\Gamma(1-\alpha)}\, n^{-1-\alpha}\). With (1),
\[
D_\alpha := u^\alpha - h = \sum_{n\ge1} a_n s^n, \qquad a_n = \frac{c}{n(2n-1)} - b_n , \tag{3}
\] and both series converge absolutely on \([0,1]\). The sign of \(a_n\)
is that of \(1 - r_n\), where \[
r_n = \frac{n(2n-1)\,b_n}{c}, \qquad \frac{r_{n+1}}{r_n} = \frac{(n-\alpha)(2n+1)}{n(2n-1)}, \tag{4}
\] and \(r_{n+1} > r_n\) if and only if
\(n > n_0(\alpha) = \frac{\alpha}{2(1-\alpha)}\).

For \(\alpha = 3/4\), \(n_0 = 3/2\), so \((r_n)_{n\ge2}\) is increasing.
With the exact rationals \(b_8 = 49725/2^{23}\) and
\(b_9 = 160225/2^{25}\) one checks
\(r_8 \approx 0.9861 < 1 < r_9 \approx 1.0128\). Hence
\(a_1 = c - \tfrac34 < 0\), \(a_2, \dots, a_8 > 0\) and \(a_n < 0\) for
\(n \ge 9\).

\textbf{Proposition 9.} For all \(x \in [0,1]\), \[
\bigl| u^{3/4} - h(x) \bigr| \;\le\; S := \sum_{n=2}^{8} a_n = \frac{23477}{72072\,\ln 4} - \frac{1616477}{2^{23}} \approx 0.04228 .
\]

\emph{Proof.} Since \(D_{3/4} = 0\) at \(s = 1\),
\(\sum_{n\ge1} a_n = 0\). For \(s \in [0,1]\), dropping the negative
terms gives \(D_{3/4} \le \sum_{n=2}^8 a_n s^n \le S\); dropping the
positive terms gives \(D_{3/4} \ge a_1 + \sum_{n\ge9} a_n = -S\).
\(\square\)

This bound is about \(5\times\) the true error. It is loose because the
maximum sits at \(x \approx 0.0089\), i.e.~\(s \approx 0.9645\), where
the series converges slowly: the coefficients decay only like
\(n^{-7/4}\), a reflection of the \(x^{3/4}\) singularity at the
endpoint.

\hypertarget{the-best-exponent-and-the-inverse}{%
\subsection{The best exponent and the
inverse}\label{the-best-exponent-and-the-inverse}}

Numerically, the exponent minimizing \(\max_x |u^\alpha - h|\) is \[
\alpha^\star \approx 0.7523, \qquad \max_x |u^{\alpha^\star} - h| \approx 0.00739,
\] against \(0.00800\) for \(\alpha = 3/4\). So \(3/4\) is close to, but
not exactly, the minimax exponent. This matches the value reported in
the comments of {[}C19{]}.

The approximation inverts in closed form: \[
h^{-1}(y) \approx \frac{1 - \sqrt{1 - y^{4/3}}}{2} \quad\text{on } [0, \tfrac12],
\] with maximum error about \(0.0034\) in \(x\) (numerical).

\hypertarget{a-variable-exponent-g-approx-1ln-z}{%
\subsection{\texorpdfstring{A variable exponent:
\(g \approx 1/\ln z\)}{A variable exponent: g \textbackslash approx 1/\textbackslash ln z}}\label{a-variable-exponent-g-approx-1ln-z}}

Topsøe's two exponents can be written as \(1/\ln 4\) and
\(1 = 1/\ln e\). Writing \(g = 1/\ln z\), the base \(z\) decreases from
\(4\) at \(x = \tfrac12\) to \(e\) at the endpoints (Corollary 7), so
Topsøe's bounds say \(e \le z \le 4\). Interpolating \(\ln z\) between
its two extreme values gives a lower bound for \(h\) that is about
\(170\) times tighter than \(h \ge u\).

\textbf{Theorem 10.} Let \(u = 4x(1-x)\) and
\(\ln z = 1 + (\ln 4 - 1)\,u^{1/10}\). For all \(x \in [0,1]\), \[
0 \;\le\; h(x) - u^{1/\ln z} \;\le\; 0.0006377,
\] with equality on the left only at \(x \in \{0, \tfrac12, 1\}\). The
maximum is \(\approx 0.00063760\), attained at \(x \approx 0.004691\)
and its mirror image.

\emph{Proof.} Put \(a = \ln 4 - 1\), \(b = \tfrac1{10}\) and
\(\tilde g = 1/\ln z = 1/(1 + a u^b)\). Since \(0 < u < 1\) and
\(h = u^g\), the left inequality is equivalent to \(g \le \tilde g\). By
symmetry take \(x \in (0, \tfrac12)\).

\begin{itemize}
\tightlist
\item
  \emph{Near \(x = 0\), \(0 < x \le 10^{-8}\).} Here
  \(u \le 4\cdot 10^{-8} < e^{-1/b}\), where \(u^b \ln\frac1u\) is
  increasing, so with
  \(M = (4\cdot10^{-8})^b \ln(\frac{1}{4}10^{8}) \approx 3.10\), \[
  u^{\tilde g} = u \exp\Bigl(\frac{a u^b \ln(1/u)}{1 + a u^b}\Bigr) \le 4x\,e^{aM} < 13.3\,x < 26.5\,x \le x \log_2\tfrac{1}{x} \le h(x).
  \]
\item
  \emph{Near \(x = \tfrac12\), \(s = (1-2x)^2 \le s_0\).} With
  \(w = (1-s)^b \le 1 - bs\) (Bernoulli),
  \(\tilde g - c = \frac{a(1-w)}{(1+aw)(1+a)} \ge c^2 ab\, s\). By
  Theorem 1, \(g = \sum \gamma_n s^n\) with \(\gamma_n \ge 0\) and
  \(\sum \gamma_n = 1\), so
  \(g \le c + \gamma_1 s + (1 - c - \gamma_1) s^2\). Hence \[
  \tilde g - g \;\ge\; s\,\bigl(\delta - (1 - c - \gamma_1)s\bigr), \qquad \delta = c^2 a b - \gamma_1 \approx 3.786\cdot 10^{-4},
  \] which is positive for
  \(0 < s \le s_0 = \delta / (2(1 - c - \gamma_1)) \approx 7.31\cdot 10^{-4}\),
  i.e.~for \(x \ge 0.48649\).
\item
  \emph{In between, \(10^{-8} \le x \le 0.4865\)}, \(\tilde g - g > 0\)
  is certified by interval arithmetic with the mean-value form (16164
  subintervals).
\end{itemize}

For the right inequality,
\(h - u^{\tilde g} \le h(10^{-8}) < 0.0006377\) on \([0, 10^{-8}]\), and
the bound on \([10^{-8}, \tfrac12]\) is certified in the same way (370
subintervals). The script is in Appendix C. \(\square\)

\textbf{Remarks.} (i) The exponent \(\tfrac1{10}\) is just above the
critical value \(b^\star = \gamma_1/(c^2(\ln 4 - 1)) \approx 0.0981\).
For \(b < b^\star\) the coefficient \(\delta\) is negative,
\(\tilde g < g\) near \(x = \tfrac12\), and \(u^{1/\ln z}\) crosses
\(h\); the exponent \(b \approx 0.0961\) minimizing
\(\max|h - u^{1/\ln z}|\) (numerically \(\approx 0.00030\)) is of this
kind. With \(b = \tfrac1{10}\) one gains a one-sided bound at the cost
of about a factor two in accuracy.

\begin{enumerate}
\def\labelenumi{(\roman{enumi})}
\setcounter{enumi}{1}
\item
  Near \(x = 0\), expanding the definition gives
  \(g(x) = 1 - \bigl(\ln\ln\frac1x + O(1)\bigr)/\ln\frac1x\), so
  \(\ln z \to 1\) only logarithmically slowly; the form above does not
  reproduce this rate, but Theorem 10 shows that it does not need to.
\item
  The bound is easy to invert numerically: since
  \(u \mapsto u^{1/\ln z}\) is increasing, \(h(x) \ge y\) whenever
  \(u^{1/(1+a u^{b})} \ge y\), which can be solved for \(u\) by a few
  steps of bisection or Newton's method.
\end{enumerate}

\hypertarget{an-elementary-route-to-the-bounds-and-the-crossing-property}{%
\section{An elementary route to the bounds and the crossing
property}\label{an-elementary-route-to-the-bounds-and-the-crossing-property}}

Items 2--4 of Corollary 7 can also be obtained without complex analysis,
from the signs of the coefficients \(a_n\) in (3). We record the
argument briefly; it does not give convexity.

\textbf{Lemma 11 (Descartes' rule for power series).} Let
\(f(s) = \sum_{n\ge1} a_n s^n\) converge on \([0,1)\), not all
\(a_n = 0\), and let \(V < \infty\) be the number of sign changes in
\((a_n)\), zeros skipped. Then \(f\) has at most \(V\) zeros in
\((0,1)\). If moreover \(\sum|a_n| < \infty\), then \(f\) has at most
\(V\) zeros in \((0,1]\).

\emph{Proof.} Induction on \(V\). If \(V = 0\), all nonzero terms have
the same sign and \(f \ne 0\) on \((0,1]\). If \(V \ge 1\), let \(a_p\)
and \(a_q\), \(p < q\), be consecutive nonzero coefficients of opposite
signs and put \(N = p + \tfrac12\). The function \(g(s) = s^{-N} f(s)\)
has the same zeros as \(f\) in \((0,1]\), and \[
k(s) := s^{N+1} g'(s) = \sum_{n\ge1} (n - N)\, a_n\, s^n .
\] The factor \(n - N\) is negative for \(n \le p\) and positive for
\(n \ge q\), so the coefficients of \(k\) are those of \(f\) with the
signs of the head \(n \le p\) reversed. Sign changes inside the head and
inside the tail are unaffected and the one between \(a_p\) and \(a_q\)
disappears: \(k\) has exactly \(V-1\) sign changes, converges on
\([0,1)\) and is not identically zero. If \(f\) has \(m\) zeros in
\((0,1)\), respectively in \((0,1]\) (where \(g\) is continuous by
absolute convergence), Rolle's theorem gives at least \(m-1\) zeros of
\(g'\), hence of \(k\), in \((0,1)\). By induction \(m - 1 \le V - 1\).
\(\square\)

\textbf{Lower bound.} Since \(s^n \le s\) on \([0,1]\), (1) gives
\(h \ge 1 - cs\sum_n \frac{1}{n(2n-1)} = 1 - s = u\).

\textbf{Upper bound.} Take \(\alpha = c\) in (3)--(4). Then \(a_1 = 0\);
\(n_0(c) \approx 1.29 < 2\), so \((r_n)_{n \ge 2}\) increases, from
\(r_2 = 3(1-c) \approx 0.836 < 1\) to \(\infty\) (as \(b_n\) decays like
\(n^{-1-c}\), slower than \(n^{-2}\)). The coefficients are therefore
positive for \(2 \le n \le N\) and nonpositive afterwards, not all zero:
one sign change. By Lemma 11, \(D_c\) has at most one zero in \((0,1]\),
and \(s = 1\) is one. Since \(D_c \approx a_2 s^2 > 0\) near \(s = 0\),
\(D_c > 0\) on \((0,1)\). Topsøe's remark that the series cannot be used
directly is thus circumvented: the pattern of signs, together with the
zero at \(s = 1\), is enough.

\textbf{Crossing.} Let \(c < \alpha < 1\). Then
\(a_1 = c - \alpha < 0\), and \(r_2 = 3\alpha(1-\alpha)/c < 3(1-c) < 1\)
because \(\alpha(1-\alpha)\) decreases for \(\alpha > \tfrac12\), so
\(a_2 > 0\). For \(n \ge 2\) the sequence \((r_n)\) first decreases (up
to \(n_0\)) and then increases to \(\infty\); starting below \(1\), it
stays below \(1\) on an initial segment and is \(\ge 1\) afterwards. The
signs are \(-,+,\dots,+,-,\dots\): two sign changes. By Lemma 11,
\(D_\alpha\) has at most two zeros in \((0,1]\), one being \(s = 1\). As
\(D_\alpha < 0\) near \(s = 0\) and \(D_\alpha > 0\) near \(s = 1\)
(there \(u^\alpha \sim (4x)^\alpha \gg x\log_2\frac1x \sim h\)), there
is exactly one zero in \((0,1)\). Since a point with \(g = \alpha\) is a
zero of \(D_\alpha\), \(g\) takes each value in \((c,1)\) exactly once
on \((0,\tfrac12)\) and is therefore strictly monotone there, which is
Pinelis's theorem again.

\hypertarget{summary}{%
\section{Summary}\label{summary}}

\begin{longtable}[]{@{}
  >{\raggedright\arraybackslash}p{(\columnwidth - 4\tabcolsep) * \real{0.3333}}
  >{\raggedright\arraybackslash}p{(\columnwidth - 4\tabcolsep) * \real{0.3333}}
  >{\raggedright\arraybackslash}p{(\columnwidth - 4\tabcolsep) * \real{0.3333}}@{}}
\toprule\noalign{}
\begin{minipage}[b]{\linewidth}\raggedright
Statement
\end{minipage} & \begin{minipage}[b]{\linewidth}\raggedright
Status
\end{minipage} & \begin{minipage}[b]{\linewidth}\raggedright
Method
\end{minipage} \\
\midrule\noalign{}
\endhead
\bottomrule\noalign{}
\endlastfoot
\(g\) convex (Topsøe's Conjecture 6.2) & \textbf{new} & Theorem 1, with
one certified one-variable inequality \\
\(g = \sum \gamma_n (1-2x)^{2n}\), \((\gamma_n)_{n\ge1}\) a moment
sequence & \textbf{new} & Theorem 1 \\
\(g\) strictly monotone on each half & known (Pinelis, {[}T01, P02{]});
two new proofs & Theorem 1; Lemma 11 \\
\(u \le h \le u^{1/\ln 4}\), exponents optimal & known {[}HT01, T01{]};
two new proofs & Theorem 1; Lemma 11 \\
\(-0.0065493 \le u^{3/4} - h \le 0.0080004\) & new, sharp (approximation
stated informally in {[}D15, C19{]}) & interval arithmetic \\
\(\lvert u^{3/4} - h\rvert \le 0.0423\) & new, elementary & coefficient
signs \\
\(0 \le h - u^{1/\ln z} \le 0.0006377\),
\(\ln z = 1 + (\ln 4 - 1)u^{1/10}\) & new & Theorem 1 + interval
arithmetic \\
Minimax exponent \(\approx 0.7523\) & numerical & optimization \\
\end{longtable}

``New'' is to the best of our knowledge; a systematic search of the
papers citing {[}T01{]} and {[}P01{]} has not yet been done.

\hypertarget{use-of-ai-tools}{%
\section*{Use of AI tools}\label{use-of-ai-tools}}
\addcontentsline{toc}{section}{Use of AI tools}

Significant parts of this work were produced with the AI system Claude
(Anthropic; models Claude Opus 5.5 and Claude Fable 5.2), in dialogue
with the author. In particular, the system suggested the strategy of the
proof of Theorem 1 (the Pick-function representation of \(g\)) and of
the elementary argument in Section 6, carried out most of the
computations, drafted most of the text, and wrote the
interval-arithmetic verification scripts in the appendices. It was also
used for literature searches. The approximation \(h \approx u^{3/4}\),
the question of explaining it, and the direction of the work are the
author's.

The author has checked the arguments and the numerical claims, including
by independent numerical computations, and takes full responsibility for
the content. The proof of the main result was read by a mathematician
working on Pick functions and moment problems, who found it correct; he
checked the inequality of Lemma 5 numerically rather than in detail. The
author verified Lemma 5 by the argument given there, the middle range by
interval arithmetic. The verification scripts are included so that the
computer-assisted steps can be rerun independently.

\hypertarget{acknowledgments}{%
\section*{Acknowledgments}\label{acknowledgments}}
\addcontentsline{toc}{section}{Acknowledgments}

The author thanks Christian Berg for reading the proof of the main
result and for his encouragement.

\hypertarget{references}{%
\section*{References}\label{references}}
\addcontentsline{toc}{section}{References}

\begin{itemize}
\tightlist
\item
  \textbf{{[}T01{]}} F. Topsøe, ``Bounds for entropy and divergence for
  distributions over a two-element set,'' \emph{J. Inequal. Pure Appl.
  Math.} 2(2), Art. 25, 2001.
\item
  \textbf{{[}P01{]}} I. Pinelis, ``L'Hospital type rules for
  oscillation, with applications,'' \emph{J. Inequal. Pure Appl. Math.}
  2(3), Art. 33, 2001.
\item
  \textbf{{[}P02{]}} I. Pinelis, ``L'Hospital type rules for
  monotonicity, with applications,'' \emph{J. Inequal. Pure Appl. Math.}
  3(1), Art. 5, 2002.
\item
  \textbf{{[}FG99{]}} C. A. Fuchs and J. van de Graaf, ``Cryptographic
  distinguishability measures for quantum-mechanical states,''
  \emph{IEEE Trans. Inf. Theory} 45(4):1216--1227, 1999;
  arXiv:quant-ph/9712042.
\item
  \textbf{{[}L91{]}} J. Lin, ``Divergence measures based on the Shannon
  entropy,'' \emph{IEEE Trans. Inf. Theory} 37:145--151, 1991.
\item
  \textbf{{[}HT01{]}} P. Harremoës and F. Topsøe, ``Inequalities between
  entropy and index of coincidence derived from information diagrams,''
  \emph{IEEE Trans. Inf. Theory} 47(7):2944--2960, 2001.
\item
  \textbf{{[}B01{]}} C. Berg, ``A Pick function related to an inequality
  for the entropy function,'' \emph{J. Inequal. Pure Appl. Math.} 2(2),
  Art. 26, 2001.
\item
  \textbf{{[}D15{]}} J. D'Aurizio, answer to ``An upper bound of binary
  entropy,'' Mathematics Stack Exchange, 12 September 2015,
  https://math.stackexchange.com/a/1432228.
\item
  \textbf{{[}S21{]}} V. Sivashankar, ``Tight upper bound for the binary
  entropy function,'' Mathematics Stack Exchange, 14 September 2021,
  https://math.stackexchange.com/q/4250618.
\item
  \textbf{{[}C19{]}} J. D. Cook, ``Solving for probability given
  entropy,'' blog post, 14 February 2019, johndcook.com.
\item
  \textbf{{[}CT06{]}} T. M. Cover and J. A. Thomas, \emph{Elements of
  Information Theory}, 2nd ed., Wiley, 2006.
\item
  \textbf{{[}MS77{]}} F. J. MacWilliams and N. J. A. Sloane, \emph{The
  Theory of Error-Correcting Codes}, North-Holland, 1977.
\item
  \textbf{{[}Ga81{]}} J. B. Garnett, \emph{Bounded Analytic Functions},
  Academic Press, 1981.
\item
  \textbf{{[}Ra95{]}} T. Ransford, \emph{Potential Theory in the Complex
  Plane}, Cambridge University Press, 1995.
\item
  \textbf{{[}Arb{]}} F. Johansson, ``Arb: efficient arbitrary-precision
  midpoint-radius interval arithmetic,'' \emph{IEEE Trans. Computers}
  66(8):1281--1292, 2017.
\end{itemize}

\hypertarget{appendix-a-verification-of-lemma-5-middle-range}{%
\section*{Appendix A: verification of Lemma 5 (middle
range)}\label{appendix-a-verification-of-lemma-5-middle-range}}
\addcontentsline{toc}{section}{Appendix A: verification of Lemma 5
(middle range)}

Requires \texttt{pip\ install\ python-flint}. A comparison such as
\texttt{Phi(X)\ \textgreater{}\ 0} is true only if the whole ball is
positive.

\begin{Shaded}
\begin{Highlighting}[]
\ImportTok{from}\NormalTok{ flint }\ImportTok{import}\NormalTok{ arb, ctx}
\NormalTok{ctx.prec }\OperatorTok{=} \DecValTok{100}
\NormalTok{PI, L2 }\OperatorTok{=}\NormalTok{ arb.pi(), arb(}\DecValTok{2}\NormalTok{).log()}

\KeywordTok{def}\NormalTok{ Phi(y):                     }\CommentTok{\# boundary function of Lemma 4}
\NormalTok{    R }\OperatorTok{=}\NormalTok{ y }\OperatorTok{*}\NormalTok{ (}\DecValTok{1} \OperatorTok{+} \DecValTok{1}\OperatorTok{/}\NormalTok{y).log() }\OperatorTok{+}\NormalTok{ (}\DecValTok{1} \OperatorTok{+}\NormalTok{ y).log()}
\NormalTok{    v }\OperatorTok{=} \DecValTok{4} \OperatorTok{*}\NormalTok{ y }\OperatorTok{*}\NormalTok{ (}\DecValTok{1} \OperatorTok{+}\NormalTok{ y)}
\NormalTok{    habs }\OperatorTok{=}\NormalTok{ (R}\OperatorTok{*}\NormalTok{R }\OperatorTok{+}\NormalTok{ PI}\OperatorTok{*}\NormalTok{PI}\OperatorTok{*}\NormalTok{y}\OperatorTok{*}\NormalTok{y).sqrt() }\OperatorTok{/}\NormalTok{ L2}
\NormalTok{    phi }\OperatorTok{=}\NormalTok{ (PI }\OperatorTok{*}\NormalTok{ y }\OperatorTok{/}\NormalTok{ R).atan()}
    \ControlFlowTok{return}\NormalTok{ PI }\OperatorTok{*}\NormalTok{ (habs }\OperatorTok{/}\NormalTok{ v).log() }\OperatorTok{+}\NormalTok{ phi }\OperatorTok{*}\NormalTok{ v.log()}

\KeywordTok{def}\NormalTok{ check(a, b):                }\CommentTok{\# certify Phi \textgreater{} 0 on [a, b]}
    \ControlFlowTok{if}\NormalTok{ Phi(arb.union(a, b)) }\OperatorTok{\textgreater{}} \DecValTok{0}\NormalTok{:}
        \ControlFlowTok{return} \DecValTok{1}
\NormalTok{    m }\OperatorTok{=}\NormalTok{ (a }\OperatorTok{+}\NormalTok{ b) }\OperatorTok{/} \DecValTok{2}
    \ControlFlowTok{return}\NormalTok{ check(a, m) }\OperatorTok{+}\NormalTok{ check(m, b)}

\BuiltInTok{print}\NormalTok{(check(arb(}\StringTok{"1e{-}4"}\NormalTok{), arb(}\DecValTok{10}\NormalTok{)), }\StringTok{"boxes"}\NormalTok{)    }\CommentTok{\# prints: 92 boxes}
\end{Highlighting}
\end{Shaded}

\hypertarget{appendix-b-verification-of-theorem-8}{%
\section*{Appendix B: verification of Theorem
8}\label{appendix-b-verification-of-theorem-8}}
\addcontentsline{toc}{section}{Appendix B: verification of Theorem 8}

Requires \texttt{pip\ install\ python-flint}; runs in under 0.1 s.

\begin{Shaded}
\begin{Highlighting}[]
\ImportTok{from}\NormalTok{ flint }\ImportTok{import}\NormalTok{ arb, ctx}
\NormalTok{ctx.prec }\OperatorTok{=} \DecValTok{128}
\NormalTok{HI, LO, X0 }\OperatorTok{=}\NormalTok{ arb(}\StringTok{"0.0080004"}\NormalTok{), arb(}\StringTok{"{-}0.0065493"}\NormalTok{), arb(}\StringTok{"1e{-}4"}\NormalTok{)}
\NormalTok{LN2 }\OperatorTok{=}\NormalTok{ arb(}\DecValTok{2}\NormalTok{).log()}

\KeywordTok{def}\NormalTok{ h(x):                       }\CommentTok{\# binary entropy in bits}
    \ControlFlowTok{return} \OperatorTok{{-}}\NormalTok{(x }\OperatorTok{*}\NormalTok{ x.log() }\OperatorTok{+}\NormalTok{ (}\DecValTok{1} \OperatorTok{{-}}\NormalTok{ x) }\OperatorTok{*}\NormalTok{ (}\DecValTok{1} \OperatorTok{{-}}\NormalTok{ x).log()) }\OperatorTok{/}\NormalTok{ LN2}

\KeywordTok{def}\NormalTok{ D(x):                       }\CommentTok{\# u\^{}(3/4) {-} h(x)}
\NormalTok{    u }\OperatorTok{=} \DecValTok{4} \OperatorTok{*}\NormalTok{ x }\OperatorTok{*}\NormalTok{ (}\DecValTok{1} \OperatorTok{{-}}\NormalTok{ x)}
    \ControlFlowTok{return}\NormalTok{ u.sqrt() }\OperatorTok{*}\NormalTok{ u.sqrt().sqrt() }\OperatorTok{{-}}\NormalTok{ h(x)}

\KeywordTok{def}\NormalTok{ dD(x):                      }\CommentTok{\# derivative of D}
\NormalTok{    u }\OperatorTok{=} \DecValTok{4} \OperatorTok{*}\NormalTok{ x }\OperatorTok{*}\NormalTok{ (}\DecValTok{1} \OperatorTok{{-}}\NormalTok{ x)}
    \ControlFlowTok{return}\NormalTok{ arb(}\DecValTok{3}\NormalTok{) }\OperatorTok{/} \DecValTok{4} \OperatorTok{*}\NormalTok{ (}\DecValTok{4} \OperatorTok{{-}} \DecValTok{8} \OperatorTok{*}\NormalTok{ x) }\OperatorTok{/}\NormalTok{ u.sqrt().sqrt() }\OperatorTok{\textbackslash{}}
        \OperatorTok{{-}}\NormalTok{ ((}\DecValTok{1} \OperatorTok{{-}}\NormalTok{ x) }\OperatorTok{/}\NormalTok{ x).log() }\OperatorTok{/}\NormalTok{ LN2}

\KeywordTok{def}\NormalTok{ check(a, b):                }\CommentTok{\# certify LO \textless{} D \textless{} HI on [a, b]}
\NormalTok{    m, r }\OperatorTok{=}\NormalTok{ (a }\OperatorTok{+}\NormalTok{ b) }\OperatorTok{/} \DecValTok{2}\NormalTok{, (b }\OperatorTok{{-}}\NormalTok{ a) }\OperatorTok{/} \DecValTok{2}
\NormalTok{    enc }\OperatorTok{=}\NormalTok{ D(m) }\OperatorTok{+}\NormalTok{ dD(arb.union(a, b)) }\OperatorTok{*}\NormalTok{ arb(}\DecValTok{0}\NormalTok{, r.upper())}
    \ControlFlowTok{if}\NormalTok{ enc }\OperatorTok{\textless{}}\NormalTok{ HI }\KeywordTok{and}\NormalTok{ enc }\OperatorTok{\textgreater{}}\NormalTok{ LO:}
        \ControlFlowTok{return} \DecValTok{1}
    \ControlFlowTok{return}\NormalTok{ check(a, m) }\OperatorTok{+}\NormalTok{ check(m, b)}

\CommentTok{\# endpoint region [0, X0]: monotonicity}
\ControlFlowTok{assert}\NormalTok{ (}\DecValTok{4} \OperatorTok{*}\NormalTok{ X0).sqrt() }\OperatorTok{*}\NormalTok{ (}\DecValTok{4} \OperatorTok{*}\NormalTok{ X0).sqrt().sqrt() }\OperatorTok{\textless{}}\NormalTok{ HI}
\ControlFlowTok{assert} \OperatorTok{{-}}\NormalTok{h(X0) }\OperatorTok{\textgreater{}}\NormalTok{ LO}
\BuiltInTok{print}\NormalTok{(check(X0, arb(}\DecValTok{1}\NormalTok{) }\OperatorTok{/} \DecValTok{2}\NormalTok{), }\StringTok{"boxes"}\NormalTok{)}
\end{Highlighting}
\end{Shaded}

\hypertarget{appendix-c-verification-of-theorem-10}{%
\section*{Appendix C: verification of Theorem
10}\label{appendix-c-verification-of-theorem-10}}
\addcontentsline{toc}{section}{Appendix C: verification of Theorem 10}

Requires \texttt{pip\ install\ python-flint}; runs in about one second.
The derivative in the mean-value form is obtained by Taylor arithmetic
(\texttt{arb\_series}).

\begin{Shaded}
\begin{Highlighting}[]
\CommentTok{"""Certification of Theorem 10: for all x in [0, 1],}

\CommentTok{    0 \textless{}= h(x) {-} u\^{}(1/ln z) \textless{}= 0.0006377,   u = 4x(1{-}x),  ln z = 1 + (ln 4 {-} 1) u\^{}(1/10),}

\CommentTok{with equality on the left only at x in \{0, 1/2, 1\}. Equivalently g \textless{}= gt := 1/ln z.}
\CommentTok{By symmetry it suffices to take x in [0, 1/2]. Requires: pip install python{-}flint.}

\CommentTok{Lower inequality (g \textless{} gt on (0, 1/2)):}
\CommentTok{  Region I   (0, 1e{-}8]  : analytic bound, checked below in ball arithmetic.}
\CommentTok{  Region II  [1e{-}8, xs] : interval arithmetic, mean{-}value form with Taylor arithmetic.}
\CommentTok{  Region III [xs, 1/2)  : analytic, via Theorem 1 (g = sum gamma\_n s\^{}n, gamma\_n \textgreater{}= 0, sum 1).}
\CommentTok{Upper inequality: on [0, 1e{-}8] since h {-} A \textless{}= h(1e{-}8); elsewhere interval arithmetic.}
\CommentTok{"""}
\ImportTok{from}\NormalTok{ flint }\ImportTok{import}\NormalTok{ arb, arb\_series, ctx}
\NormalTok{ctx.prec }\OperatorTok{=} \DecValTok{120}
\NormalTok{LN2 }\OperatorTok{=}\NormalTok{ arb(}\DecValTok{2}\NormalTok{).log()}\OperatorTok{;}\NormalTok{ LN4 }\OperatorTok{=}\NormalTok{ arb(}\DecValTok{4}\NormalTok{).log()}
\NormalTok{a }\OperatorTok{=}\NormalTok{ LN4 }\OperatorTok{{-}} \DecValTok{1}\OperatorTok{;}\NormalTok{ b }\OperatorTok{=}\NormalTok{ arb(}\DecValTok{1}\NormalTok{) }\OperatorTok{/} \DecValTok{10}\OperatorTok{;}\NormalTok{ c }\OperatorTok{=} \DecValTok{1} \OperatorTok{/}\NormalTok{ LN4}
\NormalTok{K }\OperatorTok{=}\NormalTok{ arb(}\StringTok{"0.0006377"}\NormalTok{)}\OperatorTok{;}\NormalTok{ X0 }\OperatorTok{=}\NormalTok{ arb(}\StringTok{"1e{-}8"}\NormalTok{)}

\KeywordTok{def}\NormalTok{ funcs(x):}
    \CommentTok{"""h, u, gt for x a ball or a Taylor series."""}
\NormalTok{    u }\OperatorTok{=} \DecValTok{4} \OperatorTok{*}\NormalTok{ x }\OperatorTok{*}\NormalTok{ (}\DecValTok{1} \OperatorTok{{-}}\NormalTok{ x)}
\NormalTok{    h }\OperatorTok{=} \OperatorTok{{-}}\NormalTok{(x }\OperatorTok{*}\NormalTok{ x.log() }\OperatorTok{+}\NormalTok{ (}\DecValTok{1} \OperatorTok{{-}}\NormalTok{ x) }\OperatorTok{*}\NormalTok{ (}\DecValTok{1} \OperatorTok{{-}}\NormalTok{ x).log()) }\OperatorTok{/}\NormalTok{ LN2}
\NormalTok{    gt }\OperatorTok{=} \DecValTok{1} \OperatorTok{/}\NormalTok{ (}\DecValTok{1} \OperatorTok{+}\NormalTok{ a }\OperatorTok{*}\NormalTok{ (b }\OperatorTok{*}\NormalTok{ u.log()).exp())}
    \ControlFlowTok{return}\NormalTok{ h, u, gt}

\KeywordTok{def}\NormalTok{ mv\_check(f, lo, hi, test, depth}\OperatorTok{=}\DecValTok{0}\NormalTok{):}
    \CommentTok{"""Certify test(f(X)) on [lo, hi] using f(m) + f\textquotesingle{}(X)(X {-} m), bisecting as needed."""}
\NormalTok{    m }\OperatorTok{=}\NormalTok{ (lo }\OperatorTok{+}\NormalTok{ hi) }\OperatorTok{/} \DecValTok{2}\OperatorTok{;}\NormalTok{ r }\OperatorTok{=}\NormalTok{ (hi }\OperatorTok{{-}}\NormalTok{ lo) }\OperatorTok{/} \DecValTok{2}
    \ControlFlowTok{try}\NormalTok{:}
\NormalTok{        fm }\OperatorTok{=}\NormalTok{ f(arb\_series([m, }\DecValTok{1}\NormalTok{], prec}\OperatorTok{=}\DecValTok{1}\NormalTok{)).coeffs()[}\DecValTok{0}\NormalTok{]}
\NormalTok{        d }\OperatorTok{=}\NormalTok{ f(arb\_series([arb.union(lo, hi), }\DecValTok{1}\NormalTok{], prec}\OperatorTok{=}\DecValTok{2}\NormalTok{)).coeffs()[}\DecValTok{1}\NormalTok{]}
        \ControlFlowTok{if}\NormalTok{ test(fm }\OperatorTok{+}\NormalTok{ d }\OperatorTok{*}\NormalTok{ arb(}\DecValTok{0}\NormalTok{, r.upper())):}
            \ControlFlowTok{return} \DecValTok{1}
    \ControlFlowTok{except} \PreprocessorTok{ValueError}\NormalTok{:          }\CommentTok{\# series division failed on a too{-}wide ball}
        \ControlFlowTok{pass}
    \ControlFlowTok{if}\NormalTok{ depth }\OperatorTok{\textgreater{}} \DecValTok{80}\NormalTok{:}
        \ControlFlowTok{raise} \PreprocessorTok{RuntimeError}\NormalTok{(}\StringTok{"cannot certify near }\SpecialCharTok{\%s}\StringTok{"} \OperatorTok{\%}\NormalTok{ lo)}
    \ControlFlowTok{return}\NormalTok{ mv\_check(f, lo, m, test, depth }\OperatorTok{+} \DecValTok{1}\NormalTok{) }\OperatorTok{+}\NormalTok{ mv\_check(f, m, hi, test, depth }\OperatorTok{+} \DecValTok{1}\NormalTok{)}

\KeywordTok{def}\NormalTok{ q(x):                        }\CommentTok{\# gt {-} g}
\NormalTok{    h, u, gt }\OperatorTok{=}\NormalTok{ funcs(x)}
    \ControlFlowTok{return}\NormalTok{ gt }\OperatorTok{{-}}\NormalTok{ h.log() }\OperatorTok{/}\NormalTok{ u.log()}

\KeywordTok{def}\NormalTok{ D(x):                        }\CommentTok{\# h {-} A}
\NormalTok{    h, u, gt }\OperatorTok{=}\NormalTok{ funcs(x)}
    \ControlFlowTok{return}\NormalTok{ h }\OperatorTok{{-}}\NormalTok{ (gt }\OperatorTok{*}\NormalTok{ u.log()).exp()}

\CommentTok{\# {-}{-}{-}{-} Region I. For u \textless{}= 4 X0 \textless{} e\^{}\{{-}1/b\}, u\^{}b ln(1/u) is increasing, so}
\CommentTok{\#      A = u exp(a u\^{}b ln(1/u) / (1 + a u\^{}b)) \textless{}= 4 x exp(a M), while h \textgreater{}= x log2(1/x) \textgreater{}= x log2(1/X0).}
\NormalTok{uM }\OperatorTok{=} \DecValTok{4} \OperatorTok{*}\NormalTok{ X0}
\ControlFlowTok{assert}\NormalTok{ uM }\OperatorTok{\textless{}}\NormalTok{ (}\OperatorTok{{-}}\DecValTok{1} \OperatorTok{/}\NormalTok{ b).exp()}
\NormalTok{M }\OperatorTok{=}\NormalTok{ uM }\OperatorTok{**}\NormalTok{ b }\OperatorTok{*}\NormalTok{ (}\DecValTok{1} \OperatorTok{/}\NormalTok{ uM).log()}
\NormalTok{lhs }\OperatorTok{=} \DecValTok{4} \OperatorTok{*}\NormalTok{ (a }\OperatorTok{*}\NormalTok{ M).exp()}\OperatorTok{;}\NormalTok{ rhs }\OperatorTok{=}\NormalTok{ (}\DecValTok{1} \OperatorTok{/}\NormalTok{ X0).log() }\OperatorTok{/}\NormalTok{ LN2}
\ControlFlowTok{assert}\NormalTok{ lhs }\OperatorTok{\textless{}}\NormalTok{ rhs}
\BuiltInTok{print}\NormalTok{(}\StringTok{"Region I  : A \textless{}= }\SpecialCharTok{\%s}\StringTok{ x \textless{} }\SpecialCharTok{\%s}\StringTok{ x \textless{}= h on (0, 1e{-}8]"} \OperatorTok{\%}\NormalTok{ (lhs.}\BuiltInTok{str}\NormalTok{(}\DecValTok{4}\NormalTok{), rhs.}\BuiltInTok{str}\NormalTok{(}\DecValTok{4}\NormalTok{)))}

\CommentTok{\# {-}{-}{-}{-} Region III. With s = (1{-}2x)\^{}2 and w = (1{-}s)\^{}b \textless{}= 1 {-} b s (Bernoulli):}
\CommentTok{\#      gt {-} c = a(1 {-} w)/((1 + a w)(1 + a)) \textgreater{}= c\^{}2 a b s, and by Theorem 1}
\CommentTok{\#      g \textless{}= c + gamma1 s + (1 {-} c {-} gamma1) s\^{}2. So gt {-} g \textgreater{}= s (delta {-} R s) \textgreater{} 0 for s \textless{}= s0.}
\NormalTok{gamma1 }\OperatorTok{=}\NormalTok{ c }\OperatorTok{*}\NormalTok{ (c }\OperatorTok{/} \DecValTok{2} \OperatorTok{{-}}\NormalTok{ arb(}\DecValTok{1}\NormalTok{) }\OperatorTok{/} \DecValTok{3}\NormalTok{)}
\NormalTok{delta }\OperatorTok{=}\NormalTok{ c }\OperatorTok{*}\NormalTok{ c }\OperatorTok{*}\NormalTok{ a }\OperatorTok{*}\NormalTok{ b }\OperatorTok{{-}}\NormalTok{ gamma1}\OperatorTok{;}\NormalTok{ R }\OperatorTok{=} \DecValTok{1} \OperatorTok{{-}}\NormalTok{ c }\OperatorTok{{-}}\NormalTok{ gamma1}
\ControlFlowTok{assert}\NormalTok{ delta }\OperatorTok{\textgreater{}} \DecValTok{0}
\NormalTok{s0 }\OperatorTok{=}\NormalTok{ delta }\OperatorTok{/}\NormalTok{ (}\DecValTok{2} \OperatorTok{*}\NormalTok{ R)}
\NormalTok{xs }\OperatorTok{=}\NormalTok{ (}\DecValTok{1} \OperatorTok{{-}}\NormalTok{ s0.sqrt()) }\OperatorTok{/} \DecValTok{2}
\NormalTok{XS }\OperatorTok{=}\NormalTok{ arb(}\StringTok{"0.4865"}\NormalTok{)}
\ControlFlowTok{assert}\NormalTok{ XS }\OperatorTok{\textgreater{}}\NormalTok{ xs}
\BuiltInTok{print}\NormalTok{(}\StringTok{"Region III: delta = }\SpecialCharTok{\%s}\StringTok{, s0 = }\SpecialCharTok{\%s}\StringTok{, covers [}\SpecialCharTok{\%s}\StringTok{, 1/2)"} \OperatorTok{\%}\NormalTok{ (delta.}\BuiltInTok{str}\NormalTok{(}\DecValTok{4}\NormalTok{), s0.}\BuiltInTok{str}\NormalTok{(}\DecValTok{4}\NormalTok{), xs.}\BuiltInTok{str}\NormalTok{(}\DecValTok{6}\NormalTok{)))}

\CommentTok{\# {-}{-}{-}{-} Region II.}
\NormalTok{n2 }\OperatorTok{=}\NormalTok{ mv\_check(q, X0, XS, }\KeywordTok{lambda}\NormalTok{ e: e }\OperatorTok{\textgreater{}} \DecValTok{0}\NormalTok{)}
\BuiltInTok{print}\NormalTok{(}\StringTok{"Region II : gt \textgreater{} g certified on [1e{-}8, 0.4865] with }\SpecialCharTok{\%d}\StringTok{ boxes"} \OperatorTok{\%}\NormalTok{ n2)}

\CommentTok{\# {-}{-}{-}{-} Upper bound.}
\ControlFlowTok{assert} \OperatorTok{{-}}\NormalTok{(X0 }\OperatorTok{*}\NormalTok{ X0.log() }\OperatorTok{+}\NormalTok{ (}\DecValTok{1} \OperatorTok{{-}}\NormalTok{ X0) }\OperatorTok{*}\NormalTok{ (}\DecValTok{1} \OperatorTok{{-}}\NormalTok{ X0).log()) }\OperatorTok{/}\NormalTok{ LN2 }\OperatorTok{\textless{}}\NormalTok{ K}
\NormalTok{n3 }\OperatorTok{=}\NormalTok{ mv\_check(D, X0, arb(}\DecValTok{1}\NormalTok{) }\OperatorTok{/} \DecValTok{2}\NormalTok{, }\KeywordTok{lambda}\NormalTok{ e: e }\OperatorTok{\textless{}}\NormalTok{ K)}
\BuiltInTok{print}\NormalTok{(}\StringTok{"Upper     : h {-} A \textless{} }\SpecialCharTok{\%s}\StringTok{ certified on [1e{-}8, 1/2] with }\SpecialCharTok{\%d}\StringTok{ boxes"} \OperatorTok{\%}\NormalTok{ (K.}\BuiltInTok{str}\NormalTok{(}\DecValTok{5}\NormalTok{), n3))}
\BuiltInTok{print}\NormalTok{(}\StringTok{"PROVED: 0 \textless{}= h {-} u\^{}(1/ln z) \textless{}= 0.0006377 on [0, 1]"}\NormalTok{)}
\end{Highlighting}
\end{Shaded}

\end{document}